# On some automatic continuity theorems

M.El Azhari

**Abstract.** We give characterizations of unital uniform topological algebras and saturated locally multiplicatively convex algebras by means of multiplicative linear functionals. Some automatic continuity theorems in advertibly complete uniform topological algebras are extended to a larger class of algebras. Consequences and applications are given.

**Introduction**

   A. C. Cochran introduced, in [5], the notion of saturated uniformly A-convex algebras. M. Oudadess showed, in [9], that the class of complete saturated uniformly A-convex algebras (in the sense of Cochran) is empty and introduced a new definition of saturation in uniformly A-convex algebras. In the Banach case, M. Oudadess showed that $A$ is a saturated Banach algebra iff $A$ is a unital uniform Banach algebra. A. Beddaa extended, in [1], the definition of saturation given by M. Oudadess in uniformly A-convex algebras to locally convex algebras with nonempty set of nonzero multiplicative linear functionals.
   In the second section, we give characterizations of unital uniform topological algebras and saturated locally multiplicatively convex algebras by means of multiplicative linear functionals (Theorems 2.1 and 2.2). We obtain, as a consequence, the following result: Let $(A, \|.\|)$ be a functionally continuous normed algebra with unit. Then $A$ is a saturated normed algebra iff $A$ is a uniform normed algebra (Corollary 2.4). This improves the above result of M.Oudadess. We also obtain, as a consequence, the following result of A. Beddaa: A uniform topological algebra with unit is a saturated locally multiplicatively convex algebra [1]. We show that if $A$ is a functionally continuous, saturated locally convex algebra, whose spectrum $M(A)$ is equicontinuous, then $A$ is a uniform normed algebra (Theorem 2.5). We obtain, as a consequence, the following result of A. Beddaa: The topology of a saturated locally convex Q-algebra is normable [1]. At the end of the second section, we introduce the following property: Let $(A, (p_s)_{s \in S})$ be a locally convex algebra with nonempty set $M^*(A)$ of nonzero multiplicative linear functionals, $A$ satisfies the property (P) if for all $x \in A$ and for every $s \in S$ with $p_s(x) = 1$, there exists $f_0 \in M^*(A)$ such that $|f_0(x)| \geq 1$. We show that the class of locally convex algebras satisfying property (P) contains both; the class of uniform topological algebras and the class of saturated locally convex algebras. The class of locally convex algebras satisfying property (P) is introduced to give some automatic continuity theorems.
   In the third section, we show that some automatic continuity theorems in advertibly complete uniform topological algebras, obtained by A. Beddaa, S. J. Bhatt and M. Oudadess in [2], are true if the class of advertibly complete uniform topological algebras is replaced by the larger class of locally convex algebras satisfying property (P) (Theorems 3.1, 3.2 and 3.3). Finally, as an application, we obtain that: If $(A, \|.\|)$ is a weakly regular, functionally continuous, uniform normed algebra with unit, then $\|.\|$ is the unique uniform norm on $A$ (Theorem 3.6). This improves a result of S. J. Bhatt and D. J. Karia [3, Corollary (second affirmation)].

## 1. Preliminaries

A topological algebra is an algebra over the complex field which is also a Hausdorff topological vector space such that the multiplication is separately continuous. A locally convex algebra (lc algebra) is a topological algebra whose topology is locally convex. A locally multiplicatively convex algebra [8] (lmc algebra) is a topological algebra whose topology is determined by a family $\{p_s, s \in S\}$ of submultiplicative seminorms. For each $s \in S$, let $N_s = \{x \in A, p_s(x) = 0\}$, the quotient algebra $A_s = A/N_s$ is a normed algebra in the norm $\|x_s\|_s = p_s(x)$, $x_s = x + N_s$. A uniform seminorm on an algebra $A$ is a seminorm $p$ such that $p(x^2) = p(x)^2$ for all $x \in A$. Such a seminorm is submultiplicative [6]. A uniform topological algebra (uT-algebra) is a topological algebra whose topology is determined by a family of uniform seminorms. A uniform normed algebra is a normed algebra $(A, \|.\|)$ such that $\|x^2\| = \|x\|^2$ for all $x \in A$. Let $A$ be an algebra and $x \in A$, we denote by $sp_A(x)$ the spectrum of $x$ and by $r_A(x)$ the spectral radius of $x$. A norm $\|.\|$ on an algebra $A$ is an algebra norm if $(A, \|.\|)$ is a normed algebra. An algebra norm $\|.\|$ on $A$ is semisimple [2] if the completion of $(A, \|.\|)$ is semisimple. A topological algebra is a Q-algebra [8] if the set of quasi-invertible elements is open. For a topological algebra $A$, $M(A)$ denotes the set of nonzero continuous multiplicative linear functionals on $A$. A topological algebra $A$ is strongly semisimple if for every $x \in A$, $x \neq 0$, there exists $f \in M(A)$ such that $f(x) \neq 0$. A topological algebra is weakly regular [2] if given a closed subset $F \subset M(A)$, $F \neq M(A)$, there exists a nonzero $x \in A$ such that $f(x) = 0$ for all $f \in F$. A topological algebra $A$ is functionally continuous if $M^*(A) = M(A)$. A topological algebra is weakly $\sigma^*$-compact-regular [2] if given a compact subset $K \subset M^*(A)$, $K \neq M^*(A)$, there exists a nonzero $x \in A$ such that $f(x) = 0$ for all $f \in K$. A topological algebra is advertibly complete if a Cauchy net $(x_r)_r$ in $A$ converges in $A$ whenever for some $y \in A$, $x_r + y + x_r y$ and $x_r + y + y x_r$ both converge to 0. Let $A$ be an lc algebra with unit $e$ such that $M^*(A) \neq \emptyset$, $A$ is a saturated lc algebra [1] if the topology of $A$ is determined by a family $\{p_s, s \in S\}$ of seminorms such that (1) $p_s(e) = 1$ for all $s \in S$, (2) for all $x \in A$ and $s \in S$ with $p_s(x) = 1$, there exist $f_0$, $f \in M^*(A)$ such that $|f_0(x)| = \sup\{|f(y)|, p_s(y) \leq 1\}$. Further if $p_s$ is submultiplicative for all $s \in S$, we say that $A$ is a saturated lmc algebra. A saturated normed algebra is a normed algebra $(A, \|.\|)$ with unit $e$ such that (1) $\|e\| = 1$, (2) for all $x \in A$ with $\|x\| = 1$, there exist $f_0, f \in M^*(A)$ such that $|f_0(x)| = \sup\{|f(y)|, \|y\| \leq 1\}$.

## 2. Uniform topological algebras and saturated locally convex algebras

Let $(A, (p_s)_{s \in S})$ be an lmc algebra. For each $s \in S$, let $M_s(A) = \{f \in M(A), |f(x)| \leq p_s(x)$ for all $x \in A\}$ and let $\pi_s: A \to A_s$, $\pi_s(x) = x + N_s$, be the natural homomorphism from $A$ to $A_s$. Let $f \in M(A)$ with $N_s \subset Kerf$, we may define a multiplicative linear functional $f_s$ on $A_s$ by $f_s(x_s) = f(x)$, it is clear that $f_s \in M(A_s, \|.\|_s)$. For $x \in A$, let $\hat{x}: M(A) \to \mathbb{C}$, $\hat{x}(g) = g(x)$, $g \in M(A)$, denotes the Gelfand transform of $x$.

**Theorem 2.1.** Let $(A, T)$ be a topological algebra with unit $e$. The following assertions are equivalent:
(1) $(A, T)$ is a uT-algebra;
(2) The topology $T$ is defined by a family $\{p_s, s \in S\}$ of submultiplicative seminorms such that
   (i) $p_s(e) = 1$ for all $s \in S$;
   (ii) for all $x \in A$ and $s \in S$ with $p_s(x) = 1$, there exists $f_0 \in M_s(A)$ such that $|f_0(x)| = 1$.

**Proof.** (1) $\Rightarrow$ (2): By [7, Theorem VIII.5.1], the topology $T$ of $A$ can be defined by the family $\{q_U, U \in E\}$ of seminorms, $q_U(x) = \sup\{|\hat{x}(f)|, f \in U\}$, $E$ is the set of all compact equicontinuous subsets of $M(A)$. We have $q_U(e) = 1$ for all $U \in E$. Suppose that $q_U(x) = 1$. Since the map $|\hat{x}|: M(A) \to R$, $|\hat{x}|(f) = |\hat{x}(f)|$, is continuous and $U$ is compact, there exists $f_0 \in U$ such that $1 = |\hat{x}(f_0)| = |f_0(x)|$. Clearly, $|f_0(y)| \leq q_U(y)$ for all $y \in A$, i.e., $f_0 \in M_U(A)$.

(2) $\Rightarrow$ (1): We have $f(xy - yx) = 0$ for all $x, y \in A$ and $f \in M^*(A)$, then $A$ is commutative by (ii). Let $s \in S$, $M_s(A)$ is homeomorphic to $M(\widetilde{A_s}, \|.\|_s)$ by [8, Proposition 7.5] (($\widetilde{A_s}, \|.\|_s$) is the completion of $(A_s, \|.\|_s)$). Since $(\widetilde{A_s}, \|.\|_s)$ is a commutative Banach algebra with unit, $M(\widetilde{A_s}, \|.\|_s)$ is a nonempty compact set. Let $q_s(x) = \sup\{|\hat{x}(f)|, f \in M_s(A)\}$ for all $x \in A$. If $p_s(x) = 1$, there exists $f_0 \in M_s(A)$ such that $|f_0(x)| = 1$, we have $1 = |f_0(x)| \leq q_s(x) \leq p_s(x) = 1$. So $q_s(x) = 1$. Let $x \in A$ with $p_s(x) \neq 0$. Since $p_s(p_s(x)^{-1} x) = 1$, it follows that $q_s(p_s(x)^{-1} x) = 1$, i.e., $p_s(x) = q_s(x)$ for all $x \in A$. We have $q_s(x^2) = q_s(x)^2$ for all $x \in A$, then $(A, T)$ is a uT-algebra.

**Theorem 2.2.** Let $(A, T)$ be a topological algebra with unit $e$. The following assertions are equivalent:
(1) $(A, T)$ is a saturated lmc algebra;
(2) The topology $T$ is defined by a family $\{p_s, s \in S\}$ of submultiplicative seminorms such that
  (i) $p_s(e) = 1$ for all $s \in S$;
  (ii) for all $x \in A$ and $s \in S$ with $p_s(x) = 1$, there exists $f_0 \in M^*(A)$ such that $|f_0(x)| = 1$.

**Proof.** (1) $\Rightarrow$ (2): The topology $T$ is defined by a family $\{p_s, s \in S\}$ of submultiplicative seminorms such that
(j) $p_s(e) = 1$ for all $s \in S$;
(jj) for all $x \in A$ and $s \in S$ with $p_s(x) = 1$, there exist $f_0, f \in M^*(A)$ such that $|f_0(x)| = \sup\{|f(y)|, p_s(y) \leq 1\}$.
It is easy to show that $|f(y)| \leq |f_0(x)| p_s(y)$ for all $y \in A$, then $f \in M(A)$ and $N_s \subset Ker f$. $f_s \in M(A_s, \|.\|_s)$, then $|f_s(z)| \leq \|z\|_s$ for all $z \in A_s$. Further, as $A_s$ is unital,
$1 = \sup\{|f_s(y_s)|, \|y_s\|_s \leq 1\} = \sup\{|f(y)|, p_s(y) \leq 1\} = |f_0(x)|$.
(2) $\Rightarrow$ (1): We have $f(xy - yx) = 0$ for all $x, y \in A$ and $f \in M^*(A)$, then $A$ is commutative by (ii). Let $s \in S$, $(A_s, \|.\|_s)$ is a commutative normed algebra with unit, then $M(A_s, \|.\|_s) \neq \emptyset$. Let $g_s \in M(A_s, \|.\|_s)$, $g_s o \pi_s \in M_s(A)$. We have $1 = \sup\{|g_s(y_s)|, \|y_s\|_s \leq 1\}$
$= \sup\{|g_s o \pi_s(y)|, p_s(y) \leq 1\}$. By (ii), for all $x \in A$ and $s \in S$ with $p_s(x) = 1$, there exists $f_0 \in M^*(A)$ such that $|f_0(x)| = 1 = \sup\{|g_s o \pi_s(y)|, p_s(y) \leq 1\}$. Then $(A, T)$ is a saturated lmc algebra.

**Corollary 2.3.** ([1, Proposition 2]). A uniform topological algebra with unit is a saturated lmc algebra.

**Corollary 2.4.** Let $(A, \|.\|)$ be a functionally continuous normed algebra with unit. Then the following assertions are equivalent:
(1) $A$ is a saturated normed algebra;
(2) $A$ is a uniform normed algebra.

**Theorem 2.5.** Let $(A, (p_s)_{s \in S})$ be a functionally continuous, saturated locally convex algebra, whose spectrum $M(A)$ is equicontinuous. Then $A$ is a uniform normed algebra.

**Proof.** Let $x \in A$ and $s \in S$ with $p_s(x) \neq 0$, there exist $f_0, f \in M^*(A) = M(A)$ such that $p_s(x)^{-1}|f_0(x)| = |f_0(p_s(x)^{-1}x)| = \sup\{|f(y)|, p_s(y) \leq 1\} \geq |f(e)| = 1$. Let $\|x\| = \sup\{|f(x)|, f \in M(A)\}$, we have $p_s(x) \leq \|x\|$ for all $x \in A$ and $s \in S$. Since $A$ is Hausdorff and $M(A)$ is equicontinuous, $\|.\|$ is a continuous uniform norm on $A$, then the topology of $A$ can be defined by the uniform norm $\|.\|$.

**Corollary 2.6.** Let $A$ be a saturated locally convex Q-algebra, then $A$ is a uniform normed algebra.

**Proof.** $A$ is functionally continuous and $M(A)$ is equicontinuous [7, Proposition II.7.1].

Let $(A, (p_s)_{s \in S})$ be an lc algebra with $M^*(A) \neq \emptyset$, we say that $A$ satisfies the property (P) if for all $x \in A$ and $s \in S$ with $p_s(x) = 1$, there exists $f_0 \in M^*(A)$ such that $|f_0(x)| \geq 1$. Let $A$ be an lc algebra satisfying property (P). Then $Rad(A) \subset \cap \{Kerf, f \in M^*(A)\} = \{0\}$. Also for $x, y \in A$, $xy - yx \in \cap \{Kerf, f \in M^*(A)\} = \{0\}$. Hence $A$ is commutative and semisimple. Let $(A, (p_s)_{s \in S})$ be a uT-algebra without unit. Denote by $0_s$ the zero map from $A_s$ to $C$ and by 0 the zero map from $A$ to $C$. Let $s \in S$ and $x \in A$, $p_s(x) = \|x_s\|_s = r_{\widetilde{A_s}}(x_s) = \sup\{|g(x_s)|, g \in M(\widetilde{A_s}) \cup \{0_s\}\} = \sup\{|f(x)|, f \in M_s(A) \cup \{0\}\}$ by [8, Proposition 7.5]. Let $x \in A$ and $s \in S$ with $p_s(x) = 1$. Since $M_s(A) \cup \{0\}$ is compact [8, Proposition 7.5], there exists $f_0 \in M_s(A)$ such that $1 = |\hat{x}(f_0)| = |f_0(x)|$. By the previous statement, Theorem 2.1 and the proof of Theorem 2.5, if $A$ is a uT-algebra or a saturated lc algebra, then $A$ satisfies the property (P).

**Remark.** There exists a saturated lc algebra which is not a uT-algebra [1, p.131].

## 3. Automatic continuity theorems

In this section, we extend some automatic continuity theorems in advertibly complete uniform topological algebras, obtained in [2], to lc algebras satisfying property (P).

**Theorem 3.1.** Let $(A, (p_s)_{s \in S})$ be a weakly regular, lc algebra with unit, satisfying (P). Let $B$ be an lmc algebra, and let $\phi : A \to B$ be a one-to-one homomorphism such that $C = \overline{Im\phi}$ (the closure of $Im\phi$) is a semisimple Q-algebra. If $A$ is functionally continuous, then $\phi^{-1}/Im\phi$ is continuous.

**Proof.** Let $s \in S$ and $y \in Im\phi$ with $p_s(\phi^{-1}(y)) \neq 0$. By (P) and the fact that $A$ is functionally continuous, there exists $f_0 \in M(A)$ such that $p_s(\phi^{-1}(y)) \leq |f_0(\phi^{-1}(y))|$. The map $\phi^* : M(C) \to M(A)$, $\phi^*(f) = f \circ \phi$, is well defined and $\phi^*(M(C)) = M(A)$ [2, Proof of Theorem 1]. Since $f_0 \in M(A)$, there exists $F_0 \in M(C)$ such that $f_0 = F_0 \circ \phi$. We have $p_s(\phi^{-1}(y)) \leq |f_0(\phi^{-1}(y))| = |F_0(y)| \leq r_C(y)$. Finally, $p_s(\phi^{-1}(y)) \leq r_C(y)$ for all $y \in Im\phi$ and $s \in S$. Since $C$ is a Q-algebra, $r_C$ is continuous at 0 [8, Proposition 13.5]. Then $\phi^{-1}/Im\phi$ is continuous.

**Theorem 3.2.** Let $(A, (p_s)_{s \in S})$ be a weakly $\sigma^*$-compact-regular, lc algebra with unit, satisfying (P). Let $B$ be an lc algebra, and let $\phi : A \to B$ be a one-to-one homomorphism such that $C = \overline{Im\phi}$ is a strongly semisimple Q-algebra. Then $\phi^{-1}/Im\phi$ is continuous.

**Proof.** Let $s \in S$ and $y \in Im\phi$ with $p_s(\phi^{-1}(y)) \neq 0$. By (P), there exists $f_0 \in M^*(A)$ such that $p_s(\phi^{-1}(y)) \leq |f_0(\phi^{-1}(y))|$. The map $\phi^{**} : M(C) \to M^*(A)$, $\phi^{**}(f) = f \circ \phi$, is well defined and $\phi^{**}(M(C)) = M^*(A)$ [2, Proof of Theorem 2]. Since $f_0 \in M^*(A)$, there exists $F_0 \in M(C)$ such that $f_0 = F_0 o \phi$. We have $p_s(\phi^{-1}(y)) \leq |f_0(\phi^{-1}(y))| = |F_0(y)| \leq r_C(y)$. Finally, $p_s(\phi^{-1}(y)) \leq r_C(y)$ for all $y \in Im\phi$ and $s \in S$. Since $C$ is a Q-algebra, $r_C$ is continuous at 0 [8, Proposition 13.5]. Then $\phi^{-1}/Im\phi$ is continuous.

**Theorem 3.3.** Let $(A, (p_s)_{s \in S})$ be an lc algebra satisfying (P). Let $B$ be an lc algebra, and let $\phi : A \to B$ be a one-to-one homomorphism such that $Im\phi$ is functionally continuous with continuous product and $C = \overline{Im\phi}$ is a Q-algebra. Then $\phi^{-1}/Im\phi$ is continuous.

**Proof.** Let $s \in S$ and $y \in Im\phi$ with $p_s(\phi^{-1}(y)) \neq 0$. By (P), there exists $f_0 \in M^*(A)$ such that $p_s(\phi^{-1}(y)) \leq |f_0(\phi^{-1}(y))|$. Since $Im\phi$ is functionally continuous with continuous product, it follows that $f_0 o \phi^{-1} \in M(Im\phi) = M(C)$. We have $p_s(\phi^{-1}(y)) \leq |f_0 o \phi^{-1}(y)| \leq r_C(y)$. Finally, $p_s(\phi^{-1}(y)) \leq r_C(y)$ for all $y \in Im\phi$ and $s \in S$. Since $C$ is a Q-algebra, $r_C$ is continuous at 0 [8, Proposition 13.5]. Then $\phi^{-1}/Im\phi$ is continuous.

**Theorem 3.4.** Let $(A, (p_s)_{s \in S})$ be a uT-algebra with unit. Let $B$ be an lc algebra, and let $\phi : A \to B$ be a one-to-one homomorphism such that $C = \overline{Im\phi}$ is a Q-algebra. Assume that at least one of the following holds:
(a) $A$ is weakly regular, $B$ is an lmc algebra, and $C$ is semisimple;
(b) $A$ is weakly $\sigma^*$-compact-regular, and $C$ is strongly semisimple;
(c) $Im\phi$ is functionally continuous with continuous product.
If $\phi$ is continuous, then $A$ is a uniform normed algebra.

**Proof.** By Theorems 3.1, 3.2 and 3.3, $\phi^{-1}/Im\phi$ is continuous, then $\phi$ is an homeomorphism from $A$ to $Im\phi$, and so $Im\phi$ is a uT-algebra. $C = \overline{Im\phi}$ is a uT- Q-algebra. By Corollaries 2.3 and 2.6, $C$ is a uniform normed algebra, then $A = \phi^{-1}(Im\phi)$ is a uniform normed algebra.

**Theorem 3.5.** Let $A$ be an lc Q-algebra, $B$ be a commutative complete lmc algebra. Let $\phi : A \to B$ be a one-to-one homomorphism such that $\phi^{-1}/Im\phi$ is continuous. Then $C = \overline{Im\phi}$ is a Q-algebra.

**Proof.** Since $A$ is an lc Q-algebra, $r_A$ is continuous at 0 [8, Proposition 13.5], then there exists a continuous seminorm $p$ on $A$ such that $r_{Im\phi}(y) = r_A(\phi^{-1}(y)) \leq p(\phi^{-1}(y))$ for all $y \in Im\phi$. Let $q_1 = p o \phi^{-1}$, $q_1$ is a continuous seminorm on $Im\phi$, then $q_1$ extends as a continuous seminorm $q : C \to R_+$. We have $r_C(y) \leq r_{Im\phi}(y) = r_A(\phi^{-1}(y)) \leq q(y)$ for all $y \in Im\phi$. Let $y \in C$, $y = \lim_t y_t$ for some net $(y_t)_t$ in $Im\phi$. Let $f \in M(C)$, $|f(y)| \leq |f(y - y_t)| + |f(y_t)| \leq |f(y - y_t)| + q(y_t)$, with $\lim_t(|f(y - y_t)| + q(y_t)) = q(y)$. Since $C$ is a commutative complete lmc algebra, we have $r_C(y) = \sup\{|f(y)|, f \in M(C) \cup \{0\}\} \leq q(y)$ for all $y \in C$, then $r_C$ is continuous at 0 and so $C$ is a Q-algebra by [8, Proposition 13.5].

**Remark.** Theorem 3.5 shows that the hypothesis "$A$ is complete and m-convex" is not necessary in [2, Proposition 4].

An algebra norm $\|.\|$ on an algebra $A$ is functionally continuous if $(A, \|.\|)$ is functionally continuous.

**Theorem 3.6.** Let $(A, \|.\|)$ be a weakly regular, functionally continuous, uniform normed algebra with unit. Let $|.|$ be an algebra norm on $A$.

(1) If $|.|$ is functionally continuous, then $\|.\| \leq |.|$.
(2) If $|.|$ is semisimple, then $|.|$ is functionally continuous ( $\|.\| \leq |.|$ by (1)). Further if $|.|$ is continuous, then $|.|$ is equivalent to $\|.\|$.
(3) If $|.|$ is uniform, then $|.| = \|.\|$.

**Proof.** (1) We have $M^*(A) = M(A, \|.\|) = M(A, |.|)$. Let $x \in A$,
$\|x\| = r_B(x) = \sup\{|f(x)|, f \in M(B)\}$ ($B$ is the completion of $(A, \|.\|)$)
$= \sup\{|f(x)|, f \in M(A, \|.\|)\}$
$= \sup\{|f(x)|, f \in M(A, |.|)\}$
$\leq |x|$.

(2) By Theorem 3.1, the identity map $I : (A, |.|) \to (A, \|.\|)$, $I(x) = x$, is continuous. Since $(A, \|.\|)$ is functionally continuous and $I$ is continuous, $M^*(A) = M(A, \|.\|) = M(A, |.|)$, then $|.|$ is functionally continuous.

(3) Since $(A, |.|)$ is a uniform normed algebra, the completion of $(A, |.|)$ is a uniform Banach algebra, then it is semisimple and so $|.|$ is semisimple. By (2), $|.|$ is functionally continuous, then $\|.\| \leq |.|$ by (1). $(A, |.|)$ is a weakly regular, functionally continuous, uniform normed algebra with unit. Since $\|.\|$ is uniform, it follows that $|.| \leq \|.\|$.

**Remark.** Theorem 3.6 is an improvement of [2, Corollary 2(i)] and [3, Corollary (second affirmation)].

Ecole Normale Supérieure
Avenue Oued Akreuch
Takaddoum, BP 5118, Rabat
Morocco

E-mail : mohammed.elazhari@yahoo.fr